\documentclass[11pt]{article}
\usepackage{amsfonts}
\usepackage{amsmath}
\usepackage{amssymb}
\begin{document}

\baselineskip 18pt
\title{\Large\bf On the conditional bounds for Siegel zeros}
\author{Chaohua Jia}
\date{}
\maketitle {\small \noindent {\bf Abstract.} Under a weakened
version of Hardy-Littlewood Conjecture on the number of
representations in Goldbach problem, J. H. Fei proved bounds for the
Siegel zeros. Recently G. Bhowmik and K. Halupczok generalized Fei's
result under a weaker conjecture. In the first version of this paper
on arXiv, we pointed out a defect in the paper of G. Bhowmik and K.
Halupczok, and assumed a new conjecture and used new discussion to
overcome this defect to recover their result. Afterwards, in the
second version of their paper, G. Bhowmik and K. Halupczok assumed
previous conjecture and used our new discussion to get a weaker
result. But they did not mention our paper so that we have to make
some explanation now.}

\vskip.3in
\noindent{\bf 1. Introduction}

Let $\chi$ be a Dirichlet character ${\rm mod}\,q$. It is well-known
that there is an absolute constant $c_1>0$ such that the region
$$
1-{c_1\over \log(q(|t|+2))}<\sigma
$$
contains no zero of Dirichlet $L-$function $L(s,\chi)$ unless $\chi$
is a real character, in which case $L(s,\chi)$ has at most one real
simple zero $\beta$, where $s=\sigma+it$. If such a zero $\beta$
exists, then we call it as exceptional zero and call the character
$\chi$ as exceptional character.

In 1935, A. Page proved that there is an effective absolute constant
$c_2>0$ such that for any real character $\chi\,{\rm mod}\,q\,(q\geq
3)$, if $\beta$ is a real zero of $L(s,\chi)$, then
$$
{c_2\over \sqrt{q}\log^2 q}\leq 1-\beta.
$$
In the same year, C. L. Siegel showed that for each $\varepsilon>0$,
there is a constant $c_3(\varepsilon)>0$ such that for any real
character $\chi\,{\rm mod}\,q\,(q\geq 3)$, if $\beta$ is a real zero
of $L(s,\chi)$, then
$$
{c_3(\varepsilon)\over q^\varepsilon}\leq 1-\beta, \eqno (1)
$$
where $c_3(\varepsilon)$ is an ineffective constant. The exceptional
zero is also called Siegel zero. The study for Siegel zeros is an
important topic in the number theory. Actually we need only consider
Siegel zeros for the real primitive characters.

Some people studied the connection between Siegel zeros and Goldbach
problem. One could see [3] and [2].

Write
$$
R(n)=\sum_{\substack{p_1,\,p_2\\ n=p_1+p_2}}1, \eqno (2)
$$
where $p_1,\,p_2$ denote prime numbers. Hardy and Littlewood
conjectured that for the sufficiently large even integer $n$, we
have
$$
R(n)\sim{n\over \varphi(n)}\prod_{p\not\,|\,n}\Bigl(1-{1\over
(p-1)^2}\Bigr)\cdot{n\over \log^2n},
$$
where $\varphi(n)$ is the Euler totient function. There is a
weakened version of Hardy--Littlewood Conjecture as follows.

{\bf Conjecture 1}. \emph{There is an absolute constant $c_4>0$ such
that for the even integer $n\geq 4$, we have}
$$
R(n)\geq{c_4 n\over \log^2n}.
$$

In 2016, under Conjecture 1, J. H. Fei[3] proved the following
bounds of Siegel zeros which improve the bound (1).

{\bf Theorem 1}. \emph{Suppose that Conjecture 1 holds true. Let $q$
be a prime number with $q\equiv 3\,({\rm mod}\,4),\,\chi$ be the
real primitive character ${\rm mod}\,q$ and $\beta$ be the real zero
of $L(s,\chi)$. Then there is an absolute constant $c_5>0$ such
that}
$$
{c_5\over \log^2q}\leq 1-\beta.
$$

In Fei's device, he calculated the sum
$$
S=\sum_{k=1}^q\Bigl(\sum_{2<p\leq x}e({kp\over q})\Bigr)^2 \eqno (3)
$$
in two different ways.

In the first way, he got a lower bound for the sum in (3) by
Conjecture 1. In the second way, he connected the sum in (3) with
Siegel zero by the application of the prime number theorem in the
arithmetic progression which is Lemma 2 below. Then he obtained an
upper bound for the sum in (3) so that the bound for Siegel zero
follows.

Recently, in the first version of paper [1], G. Bhowmik and K.
Halupczok generalized the result in Theorem 1 under the following
Conjecture 2, which is shown in the following Theorem 2.

{\bf Conjecture 2}. \emph{Suppose that $x$ is sufficiently large,
$q\leq{x\over 4}$. There is an absolute constant $c_6>0$ such that
for the even integers $n({x\over 2}<n\leq x, q|\,n)$, with at most
${x\over 8q}$ exceptions, we have}
$$
R(n)\geq{c_6 n\over \log^2n}.
$$

{\bf Theorem 2}. \emph{Suppose that Conjecture 2 holds true. Let $q$
be a sufficiently large integer, $\chi$ be the real primitive
character ${\rm mod}\,q$ with $\chi(-1)=-1$ and $\beta$ be the
Siegel zero of $L(s,\chi)$.  Then there is an effective absolute
constant $c_7>0$ such that}
$$
{c_7\over \log^2q}\leq 1-\beta.
$$

G. Bhowmik and K. Halupczok worked in the framework of Fei[3]. In
the first way, they used Conjecture 2 weaker than Conjecture 1 to
replace Conjecture 1 and obtained a lower bound similar as Fei's. In
the second way, they followed Fei's discussion. But there is a
defect arising in their formula (3.6).

In Fei's discussion, $q$ is a prime number. Then for the integer
$k(1\leq k\leq q-1)$, the Ramanujan sum
$$
c_q(k)=\sum_{\substack{a=1\\ (a,\,q)=1}}^q e({ak\over q})
$$
equals to $\mu(q)$ which is bounded. But for the proper factors of
composite number $q$, $|c_q(k)|$ may be large so that the estimate
$$
O\Bigl({x\over \varphi(q)\log x}\Bigr)
$$
in (3.6) can not be arrived. We pointed out this defect in Remark 1
in the first version of this paper.

If one supposes that $|c_q(k)|$ is bounded, then ${q\over
\varphi(q)}$ would be large so that Conjecture 2 is not enough to
get Theorem 2. After a careful analysis, one could see that G.
Bhowmik and K. Halupczok actually proved Theorem 2 only for the
situation of prime numbers under a weaker conjecture.

In the first version of this paper which was published on arXiv on
October 28, 2020, we assumed a new conjecture which is a little
stronger than Conjecture 2 and used a new discussion to overcome
this defect to recover their result.

{\bf Conjecture 3}. \emph{Suppose that $x$ is sufficiently large,
$q\leq{x\over 4}$. There is an absolute constant $c_8>0$ such that
for the even integers $n({x\over 2}<n\leq x, q|\,n)$, we have}
$$
R(n)\geq{c_8n\over \varphi(n)}\cdot{n\over \log^2n}.
$$

Conjecture 3 is a little stronger than Conjecture 2, but is still a
weakened version of Hardy--Littlewood Conjecture. In the first way,
under Conjecture 3, we obtained a lower bound similar as Fei's but
with an additional factor ${q\over \varphi(q)}$. In the second way,
we used a new discussion to deal with sums involving Ramanujan sum
for composite numbers and obtained an upper bound with an additional
factor ${q\over \varphi(q)}$. Canceling the factor ${q\over
\varphi(q)}$ in both sides, we got the following Theorem 3.

{\bf Theorem 3}. \emph{Suppose that Conjecture 3 holds true. Let $q$
be a sufficiently large composite number, $\chi$ be the real
primitive character ${\rm mod}\,q$ with $\chi(-1)=-1$ and $\beta$ be
the Siegel zero of $L(s,\chi)$.  Then there is an effective absolute
constant $c_9>0$ such that}
$$
{c_9\over \log^2q}\leq 1-\beta.
$$

Some days after the publication of the first version of our paper,
G. Bhowmik and K. Halupczok published the second version of their
paper [1] on arXiv, in which they proved the following Theorem 4.

{\bf Theorem 4}. \emph{Suppose that Conjecture 2 holds true. Let $q$
be a sufficiently large integer, $\chi$ be the real primitive
character ${\rm mod}\,q$ with $\chi(-1)=-1$ and $\beta$ be the
Siegel zero of $L(s,\chi)$.  Then there is an effective absolute
constant $c_{10}>0$ such that}
$$
c_{10}{\varphi(q)\over q\log^2q}\leq 1-\beta.
$$

In the first way, they still assumed the previous Conjecture 2. In
the second way, they used our new discussion. Therefore they got a
weaker result. Since they did not mention our paper, we have to make
some explanation now.

Throughout this paper, we assume that $c_i$ is the positive
constant. Let $p,\,p_i$ denote the prime numbers, $\varphi(n)$
denote the Euler totient function, $\mu(n)$ denote the M\"obius
function.

\vskip.3in
\noindent{\bf 2. Some lemmas}

{\bf Lemma 1}. \emph{We have}
$$
{n\over \varphi(n)}\geq\sum_{d|\,n}{\mu^2(d)\over d}\geq{6\over
\pi^2}\cdot{n\over \varphi(n)}.
$$

Proof. Firstly, by the expression of the Euler totient function, we
have
\begin{align*}
{n\over \varphi(n)}&={1\over \prod_{p|\,n}(1-{1\over p})}\\
&={\prod_{p|\,n}(1+{1\over p})\over \prod_{p|\,n}(1-{1\over
p^2})}\\
&\geq\prod_{p|\,n}\Bigl(1+{1\over p}\Bigr)\\
&=\sum_{d|\,n}{\mu^2(d)\over d}.
\end{align*}

Secondly, in a similar way to the above, we get
\begin{align*}
\sum_{d|\,n}{\mu^2(d)\over d}&=\prod_{p|\,n}\Bigl(1+{1\over
p}\Bigr)\\
&={\prod_{p|\,n}(1-{1\over p^2})\over \prod_{p|\,n}(1-
{1\over p})}\\
&\geq\prod_p\Bigl(1-{1\over p^2}\Bigr)\cdot{n\over \varphi(n)}\\
&={6\over \pi^2}\cdot{n\over \varphi(n)},
\end{align*}
where the fact
$$
\prod_p\Bigl(1-{1\over p^2}\Bigr)^{-1}=\sum_{n=1}^\infty{1\over
n^2}={\pi^2\over 6}
$$
is used.

{\bf Lemma 2}. \emph{Let $\chi$ be the real primitive character
${\rm mod}\,q\,(q\geq 3),\,\beta$ be the Siegel zero of $L(s,\chi)$.
If $(a,\,q)=1$, then there is an effective absolute constant
$c_{11}>0$ such that}
$$
\sum_{\substack{p\leq x\\ p\equiv a\,({\rm mod}\,q)}}1={{\rm
li}(x)\over \varphi(q)}-{\chi(a)\over \varphi(q)}\,{\rm
li}(x^\beta)+O(x\exp(-c_{11}\sqrt{\log x})),
$$
\emph{where}
$$
{\rm li}(x)=\int_2^x{du\over \log u}.
$$

One could see Corollary 11.20 in page 381 of [4].

The sum
$$
c_q(n)=\sum_{\substack{a=1\\ (a,\,q)=1}}^q e({an\over q})
$$
is called as Ramanujan sum.

{\bf Lemma 3}. \emph{For the Ramanujan sum, we have}
$$
c_q(n)={\mu({q\over (q,\,n)})\over \varphi({q\over (q,\,n)})}\cdot
\varphi(q).
$$

One could see Theorem 4.1 in page 110 of [4].

For the character $\chi\,{\rm mod}\,q$, we write
$$
\tau(\chi,k)=\sum_{a=1}^q\chi(a)e({ak\over q})
$$
and
$$
\tau(\chi)=\tau(\chi,1).
$$

{\bf Lemma 4}. \emph{If $\chi$ is a primitive character ${\rm
mod}\,q$, then we have that}
$$
\tau(\chi,k)=\overline{\chi}(k)\tau(\chi)
$$
and
$$
|\tau(\chi)|=\sqrt{q}.
$$

One could see Theorem 9.5 and 9.7 in page 287 of [4].

{\bf Lemma 5}. \emph{If $\chi$ is the real primitive character ${\rm
mod}\,q$, then we have}
$$
\tau^2(\chi)=\chi(-1)q.
$$

Proof. It is easy to see
\begin{align*}
\overline{\tau(\chi)}&=\sum_{a=1}^q\chi(a)e(-{a\over
q})\\
&=\chi(-1)\sum_{b=1}^q\chi(b)e({b\over q})\\
&=\chi(-1)\tau(\chi).
\end{align*}
Then by Lemma 4, we have
$$
\tau^2(\chi)=\chi(-1)|\tau(\chi)|^2=\chi(-1)q.
$$

\vskip.3in
\noindent{\bf 3. The proof of Theorem 3}

Following the research route in [3], we shall consider the lower
bound and upper bound of the sum $S$ in (3).

Firstly we take
$$
x=\exp\Bigl({36\over c_{11}^2}\log^2q\Bigr),  \eqno(4)
$$
where $c_{11}$ is the constant defined in Lemma 2. Since $q$ is
sufficiently large, so is $x$. It is easy to see
$$
q=\exp\Bigl({c_{11}\over 6}\sqrt{\log x}\bigr)<{x\over 4}.
$$

We have
\begin{align*}
S&=\sum_{k=1}^q\sum_{2<p_1,\,p_2\leq x}e({k(p_1+p_2)\over q})\\
&=\sum_{n\leq 2x}\sum_{k=1}^q e({kn\over q})
\sum_{\substack{2<p_1,\,p_2\leq x\\ n=p_1+p_2}}1\\
&=\sum_{\substack{n\leq 2x\\ q|\,n\\ 2|\,n}} q
\sum_{\substack{2<p_1,\,p_2\leq x\\ n=p_1+p_2}}1
\end{align*}
\begin{align*}
&\geq q\sum_{\substack{{x\over 2}<n\leq x\\ q|\,n\\
2|\,n}}\sum_{\substack{p_1,\,p_2\\ n=p_1+p_2}}1\\
&=q\sum_{\substack{{x\over 2}<n\leq x\\ q|\,n\\ 2|\,n}}R(n).
\end{align*}

Under the assumption of Conjecture 3, for the even integers
$n({x\over 2}<n\leq x, q|\,n)$, we have
$$
R(n)\geq{c_8n\over \varphi(n)}\cdot{n\over \log^2n}.
$$
Hence, by Lemma 1,
\begin{align*}
S&\gg{qx\over \log^2x}\sum_{\substack{{x\over 2}<n\leq x\\ q|\,n\\
2|\,n}}{n\over \varphi(n)}\\
&\geq{qx\over \log^2x}\sum_{\substack{{x\over 2}<n\leq x\\ 2q|\,n}}
{n\over \varphi(n)}\\
&\geq{qx\over \log^2x}\sum_{\substack{{x\over 2}<n\leq x\\
2q|\,n}}\sum_{d|\,n}{\mu^2(d)\over d}\\
&={qx\over \log^2x}\sum_{d\leq x}{\mu^2(d)\over
d}\sum_{\substack{{x\over 2}<n\leq x\\ 2q|\,n\\ d|\,n}}1\\
&\geq{qx\over \log^2x}\sum_{\substack{d\leq x\\
d|\,q}}{\mu^2(d)\over
d}\sum_{\substack{{x\over 2}<n\leq x\\ 2q|\,n}}1\\
&\geq{qx\over \log^2x}\sum_{d|\,q}{\mu^2(d)\over d}\sum_{\substack{
{x\over 2}<n\leq x\\ 2q|\,n}}1\\
&\gg{qx\over \log^2x}\cdot{x\over 4q}\sum_{d|\,q}{\mu^2(d)\over d}\\
&\gg{q\over \varphi(q)}\cdot{x^2\over \log^2 x}.
\end{align*}
Therefore there is an absolute constant $c_{12}(0<c_{12}<{1\over
4})$ such that
$$
S\geq c_{12}\cdot{q\over \varphi(q)}\cdot{x^2\over \log^2x}. \eqno
(5)
$$

On the other hand, the application of Lemmas 2, 3 and 4 produces
that
\begin{align*}
\sum_{2<p\leq x}e({kp\over q})
&=\sum_{\substack{2<p\leq x\\
p\not\,|\,q}}e({kp\over q})+\sum_{\substack{2<p\leq x\\ p|\,q}}
e({kp\over q})\\
&=\sum_{\substack{p\leq x\\ (p,\,q)=1}}e({kp\over q})+O(\log q)\\
&=\sum_{\substack{a=1\\ (a,\,q)=1}}^q e({ka\over
q})\sum_{\substack{p\leq x\\ p\equiv a\,({\rm mod}\,q)}}1+O(\log q)\\
&=\sum_{\substack{a=1\\ (a,\,q)=1}}^q e({ka\over q})\cdot{{\rm
li}(x)\over \varphi(q)}-\sum_{\substack{a=1\\ (a,\,q)=1}}^q
\chi(a)e({ka\over q})\cdot{{\rm li}(x^\beta)\over \varphi(q)}\\
&\qquad+O(qx\exp(-c_{11}\sqrt{\log x})) \\
&=c_q(k)\cdot{{\rm li}(x)\over \varphi(q)}-\tau(\chi,k)
\cdot{{\rm li}(x^\beta)\over \varphi(q)}\\
&\qquad+O(qx\exp(-c_{11}\sqrt{\log x}))\\
&={\mu({q\over (q,\,k)})\over \varphi({q\over (q,\,k)})}\cdot{\rm
li}(x)-\chi(k)\tau(\chi)\cdot{{\rm li}(x^\beta)\over \varphi(q)}\\
&\qquad+O(qx\exp(-c_{11}\sqrt{\log x})).
\end{align*}
Hence,
\begin{align*}
&\ \,\sum_{k=1}^q\Bigl(\sum_{2<p\leq x}e({kp\over q})\Bigr)^2\\
&=\sum_{k=1}^q{\mu^2({q\over (q,\,k)})\over \varphi^2({q\over
(q,\,k)})}\cdot{\rm li}^2(x)+\sum_{k=1}^q\chi^2(k)\tau^2(\chi)
\cdot{{\rm li}^2(x^\beta)\over \varphi^2(q)}\\
&\ -2\sum_{k=1}^q{\mu({q\over (q,\,k)})\over \varphi({q\over
(q,\,k)})}\cdot{\rm li}(x)\cdot\chi(k)\tau(\chi)\cdot{{\rm
li}(x^\beta)\over \varphi(q)}\\
&\quad+O\Bigl(\sum_{k=1}^q q^2x^2\exp(-c_{11}\sqrt{\log x})\Bigr).
\end{align*}

{\bf Remark 1}. \emph{When $q$ is a prime number and $1\leq k\leq
q-1$, we have that $(k,\,q)=1$ and $c_q(k)=\mu(k)=O(1)$ as done in
[3]. But when $q$ is a composite number, we can not ensure that
$c_q(k)=O(1)$ for $1\leq k\leq q-1$. Actually for the proper factors
of $q,\,|c_q(k)|$ may be large. This is the defect in the proof of
Theorem 11 in the first version of [1].}

\emph{If one supposes that $|c_q(k)|$ is bounded, then ${q\over
\varphi(q)}$ would be large so that Conjecture 2 is not enough to
get Theorem 2.}

Note that
$$
{\rm li}(x)={x\over \log x}+O\Bigl({x\over \log^2x}\Bigr).
$$
Then we have
\begin{align*}
&\ \,\sum_{k=1}^q{\mu^2({q\over (q,\,k)})\over \varphi^2({q\over
(q,\,k)})}\cdot{\rm li}^2(x)\\
&=\sum_{d|\,q}\Bigl(\sum_{\substack{k=1\\ (q,\,k)=d}}^q 1\Bigr)
{\mu^2({q\over d})\over \varphi^2({q\over d})}\cdot{\rm li}^2(x)\\
&=\sum_{d|\,q}\Bigl(\sum_{\substack{l=1\\ (l,\,{q\over
d})=1}}^{q\over d} 1\Bigr){\mu^2({q\over d})\over \varphi^2({q\over
d})}\cdot{\rm li}^2(x)\\
&=\sum_{d|\,q}{\mu^2({q\over d})\over \varphi({q\over d})}\cdot
{\rm li}^2(x)\\
&=\sum_{r|\,q}{\mu^2(r)\over \varphi(r)}\cdot{\rm li}^2(x)\\
&=\prod_{p|\,q}\Bigl(1+{1\over \varphi(p)}\Bigr)\Bigl({x^2\over
\log^2x}+O\Bigl({x^2\over \log^3x}\Bigr)\Bigr)\\
&={q\over \varphi(q)}\cdot{x^2\over \log^2x}+O\Bigl({q\over
\varphi(q)}\cdot{x^2\over \log^3x}\Bigr),
\end{align*}
and
\begin{align*}
&\ \,\sum_{k=1}^q\chi^2(k)\tau^2(\chi)\cdot{{\rm li}^2(x^\beta)
\over \varphi^2(q)}\\
&={\chi(-1)q\over \varphi(q)}\cdot{\rm li}^2(x^\beta)\\
&={\chi(-1)q\over \varphi(q)}\Bigl({x^\beta\over \beta\log
x}+O\Bigl({x^\beta\over \log^2x}\Bigr)\Bigr)^2\\
&={\chi(-1)q\over \varphi(q)}\cdot{x^{2\beta}\over
\beta^2\log^2x}+O\Bigl({q\over \varphi(q)}\cdot{x^2\over
\log^3x}\Bigr).
\end{align*}
We also have
\begin{align*}
&-2\sum_{k=1}^q{\mu({q\over (q,\,k)})\over \varphi({q\over
(q,\,k)})}\cdot{\rm li}(x)\cdot\chi(k)\tau(\chi)\cdot{{\rm
li}(x^\beta)\over \varphi(q)}\\
&=-{2\over \varphi(q)}\cdot{\rm li}(x){\rm
li}(x^\beta)\tau(\chi)\sum_{\substack{k=1\\
(k,\,q)=1}}^q{\mu(q)\over \varphi(q)}\cdot\chi(k)\\
&=0.
\end{align*}

Combining the above estimates, we obtain
\begin{align*}
S&={q\over \varphi(q)}\cdot{x^2\over \log^2x}+{\chi(-1)q\over
\varphi(q)}\cdot{x^{2\beta}\over \beta^2\log^2x}+O\Bigl({q\over
\varphi(q)}\cdot{x^2\over \log^3x}\Bigr)\\
&\qquad\qquad\qquad+O(q^3x^2\exp(-c_{11}\sqrt{\log x}))\qquad
\qquad\qquad\qquad\qquad (6)\\
&={q\over \varphi(q)}\cdot{x^2\over \log^2x}-{q\over
\varphi(q)}\cdot{x^{2\beta}\over \beta^2\log^2x}+O\Bigl({q\over
\varphi(q)}\cdot{x^2\over \log^3x}\Bigr).
\end{align*}

Comparing (5) with (6), we get
$$
c_{12}\cdot{q\over \varphi(q)}\cdot{x^2\over \log^2x}\leq{q\over
\varphi(q)}\cdot{x^2\over \log^2x}-{q\over
\varphi(q)}\cdot{x^{2\beta}\over \beta^2\log^2x}+O\Bigl({q\over
\varphi(q)}\cdot{x^2\over \log^3x}\Bigr),
$$
which yields
$$
{x^{2\beta-2}\over \beta^2}\leq (1-c_{12})+{c_{12}\over 2}.
$$
Therefore we have
$$
x^{2\beta-2}\leq 1-{c_{12}\over 2},
$$
which yields
$$
1-\beta\geq{-\log(1-{c_{12}\over 2})\over 2\log x}
={-c_{11}^2\log(1-{c_{12}\over 2})\over 72\log^2q}={c_9\over
\log^2q},
$$
where $c_9>0$ is an effective absolute constant.

So far the proof of Theorem 3 is complete.

{\bf Remark 2}. \emph{In the proofs of Theorems 2, 3 and 4, one
needs only suppose that Conjectures 2 and 3 hold true for
$q\leq\exp(C\sqrt{\log x})$, where $C>0$ is some large constant.}

\vskip.3in
\noindent{\bf Acknowledgements}

This work is supported by the National Natural Science Foundation of
China (Grant No. 11771424).

\pagebreak

\vskip.6in

Chaohua Jia

Institute of Mathematics, Academy of Mathematics and Systems
Science, Chinese Academy of Sciences, Beijing 100190, P. R. China

Hua Loo-Keng Key Laboratory of Mathematics, Chinese Academy of
Sciences, Beijing 100190, P. R. China

School of Mathematical Sciences, University of Chinese Academy of
Sciences, Beijng 100049, P. R. China

E-mail: {\tt jiach@math.ac.cn}

\end{document}